\documentclass[11pt]{article}

\newcommand{\bc}{\begin{center}}
\newcommand{\ec}{\end{center}}
\newcommand{\be}{\begin{equation}}
\newcommand{\ee}{\end{equation}}
\newcommand{\ba}{\begin{array}}
\newcommand{\ea}{\end{array}}
\newcommand{\bea}{\begin{eqnarray}}
\newcommand{\eea}{\end{eqnarray}}
\newcommand{\edc}{\end{document}}

\def\s{\sigma}

\def\s{\sigma}

\def\l{\lambda}

\def\P{\Phi}

\def \L {\Lambda}

\begin{document}
\thispagestyle{empty}
\begin{center}
\large
{\textbf{QUADRATIC STOCHASTIC OPERATORS:\\ RESULTS AND OPEN PROBLEMS}}\\[4mm]
 {\bf R.N. Ganikhodzhaev}$^1$,  {\bf U. A. Rozikov}$^{2}$ \\[2mm]
 $^1$ Department of Mechanics and Mathematics, National University
 of Uzbekistan.\\
$^2$Institute of Mathematics and Information Technologies, Tashkent, Uzbekistan.\\
\vspace{0.3cm}
\end{center}
{\bf Abstract.} The history of the quadratic stochastic operators
can be traced back to work of S.Bernshtein (1924). During more than
80 years this theory developed and many papers were published. In
recent years it has again become of interest in connection with
numerous applications to many branches of mathematics, biology and
physics. But most results of the theory were published in non
English journals, full text of which are not accessible. In this
paper we give a brief description of the results and discuss
several open problems. \\


\section{Introduction} Quadratic stochastic operator (QSO) was first
introduced in [1].  A QSO has meaning of a population evolution
operator, which arises as follows. Consider a population consisting
of $m$ species. Let $x^{0} = (x_{1}^{0},...,x_{m}^{0})$ be the
probability distribution of species in the initial generations, and
$P_{ij,k}$ the probability that individuals in the $i$th and $j$th
species interbreed to produce an individual $k$. Then the
probability distribution $x'= (x_{1}',...,x_{m}')$ (the state) of
the species in the first generation can be found by the total
probability i.e.
$$ x'_k=\mathop {\sum} \limits^{m}_{i,j=1}P_{ij,k}x^{0}_{i}x^{0}_{j} ,\,\,\,\,k=
1,...,m. \eqno(1)$$
 This means that the association $x^{0}\rightarrow x'$
defines a map $V$ called the evolution operator. The population
evolves by starting from an arbitrary state $x^{0}$, then passing to
the state $x'= V(x)$ (in the next "generation"), then to the state
$x''=V(V(x))$, and so on.
 Thus states of the population described by the following dynamical
 system
         $$ x^{0},\ \ x'= V(x), \ \ x''=V^{2}(x),\ \  x'''= V^{3}(x),..., \eqno(2)$$
         where $V^n(x)=V(V(...V(x))...)$ denotes the $n$ times
         iteration of $V$ to $x$.

Note that $V$ (defined by (1)) is a non linear (quadratic) operator,
and it is higher dimensional if $m\geq 3$. Higher dimensional
dynamical systems are important but there are relatively few
dynamical phenomena that are currently understood [2].

The main problem for a given dynamical system (2) is to describe the
limit points of $\{x^{(n)}\}_{n=0}^\infty$ for arbitrary given
$x^{(0)}$.  In this paper we discuss recently obtained results on
the problem, also give several open problems related to theory of
QSOs.
\newpage
\section{Definitions} The quadratic stochastic operator (QSO) is a
mapping of the simplex.

$$ S^{m-1}=\left\{x=(x_1,...,x_m)\in {\bf R}^m: x_i\geq
0,\,\sum^m_{i=1}x_i=1 \right\} \eqno(3)$$ into itself , of the form
$$
    V: x'_k= \sum^m_{i,j=1}P_{ij,k}x_ix_j,\,\,\,k=
    1,...,m,\eqno(4)$$
where $P_{ij,k}$ are coefficients of heredity and
$$P_{ij,k}\geq 0,\ \ P_{ij,k}= P_{ji,k},\ \ \sum^m_{k=1}P_{ij,k}= 1, (i,j,k=
1,....,m).\eqno(5)$$

 Thus each quadratic stochastic operator $V$ can
be uniquely defined by a cubic matrix
$\textbf{P}=\left(P_{ij,k}\right)^n _{i,j,k=1}$ with conditions (5).

Note that each element $x\in S^{m-1}$ is a probability distribution
on $E= \{1,...,m\}$.

For a given $x^{(0)}\in S^{m-1}$ the trajectory (orbit)
$$ \{ x^{(n) }\}, \ \ n=0,1,2,...\ \ \mbox{of}\ \ x^{(0)}$$
under the action of QSO (4) is defined by
$$x^{(n+1)}= V(x^{(n)}),\ \ \mbox{where} \ \ n= 0,1,2,...$$

One of the main problem in mathematical biology consists in the
study of the asymptotical behavior of the trajectories. The
difficulty of the problem depends on given matrix \textbf{P}.

\section{The Volterra operators} A Volterra QSO is defined by (4),
(5) and the additional assumption
$$ P_{ij,k}=0, \ \ \mbox{if}\ \ k\not\in \{i,j\},\,\,
\forall i,j,k\in E. \eqno(6)$$

The biological treatment of condition (6) is clear: The offspring
repeats the genotype of one of its parents.

In paper [10] the general form of Volterra QSO
$$ V: x= (x_{1},...,x_{m}) \in S^{m-1}\,
\rightarrow\,V(x)= x'= (x'_1,...,x'_m)\,\in S^{m-1}
$$ is given
$$x'_k=x_k\left(1+\sum^m_{i=1}a_{ki}x_i\right),\ \ k\in E\eqno(7)$$
where
$$a_{ki}=2P_{ik,k}-1 \ \ \mbox{for}\, i\neq k\,\,\mbox{and}\
\,a_{ii}=0, \ i\in E. \eqno(8)$$  Moreover
$$a_{ki}=-a_{ik}\ \ \mbox{and} \ \ |a_{ki}| \leq 1.$$
Denote by $A=(a_{ij})_{i,j=1}^m$ the skew-symmetric matrix with
entries (8).

Let $\{x^{(n)}\}_{n=1}^\infty$ be the trajectory of the point
$x^0\in S^{m-1}$ under QSO (7). Denote by $\omega(x^0)$ the set of
limit points of the trajectory. Since $\{x^{(n)}\}\subset S^{m-1}$
and $S^{m-1}$ is compact, it follows that $\omega(x^0)\ne
\emptyset.$ Obviously, if $\omega(x^0)$ consists of a single point,
then the trajectory converges, and $\omega(x^0)$ is a fixed point of
(7). However, looking ahead, we remark that convergence of the
trajectories is not the typical case for the dynamical systems (7).
Therefore, it is of particular interest to obtain an upper bound for
$\omega(x^0)$, i.e., to determine a sufficiently "small" set
containing $\omega(x^0)$.

Denote $${\rm int}S^{m-1}=\{x\in S^{m-1}: \prod_{i=1}^m x_i>0\}, \ \
\partial S^{m-1}=S^{m-1}\setminus {\rm int}S^{m-1}.$$
 \vskip 0.2 truecm
{\bf Definition 1.} A continuous function $\varphi:{\rm
int}S^{m-1}\to R$ is called a Lyapunov function for the dynamical
system (7) if the limit $\lim_{n\to \infty}\varphi(x^{(n)})$ exists
for any initial point $x^0$. \vskip 0.2 truecm

 Obviously, if
$\lim_{n\to \infty}\varphi(x^{(n)})=c$, then $\omega(x^0)\subset
\varphi^{-1}(c).$ Consequently, for an upper estimate of
$\omega(x^0)$ we should construct the set of Lyapunov functions that
is as large as possible. \vskip 0.2 truecm

  In [3],[10]- [14], [45] the theory of QSOs (7) was developed
by using theory of the Lyapunov function and tournaments.

 The following results are known:
\vskip 0.2 truecm

 {\bf Theorem 1.} [10],[14] {\it For the Volterra
QSO (7) the following assertions hold}

1) {\it For the dynamical system (7) there exists a Lyapunov
function of the form $\varphi_p(x)=x_1^{p_1}...x_m^{p_m},$ where
$p_i\geq 0$,  $\sum^m_{i=1}p_i=1$ and  $x=(x_1,...,x_m)\in {\rm
int}S^{m-1}$.}

2) {\it If there is $r\in \{1,...,m\}$ such that $a_{ij}<0$ (see
(8)) for all $i\in \{1,...,r\}, \ j\in \{r+1,...,m\}$ then
$\varphi(x)=\sum^m_{i=r+1}x_i, \ x\in S^{m-1}$ is a Lyapunov
function for QSO (7).}

3) {\it There are Lyapunov functions of the form}
$$\varphi(x)=\frac{x_i}{x_j}, \ i\ne j, \ x\in {\rm int}S^{m-1}.$$
\vskip 0.2 truecm

 {\bf Problem 1.} {\sl Does there exist another
kind of Lyapunov function for QSO (7)?}
\vskip 0.2 truecm

The next theorem related to the set of limit points of QSO (7).
\vskip 0.2 truecm

{\bf Theorem 2.} [10],[14] 1) {\it If $x^{(0)}\in {\rm int}S^{m-1}$
is not a fixed point (i.e. $V(x^{(0)})\ne x^{(0)}$), then
$\omega(x^0)\subset \partial S^{m-1}$.}

2) {\it The set $\omega(x^0)$ either consists of a single point or
is infinite.}

3) {\it If QSO (7) has an isolated fixed point $x^*\in {\rm int}
S^{m-1}$ then for any initial point $x^{(0)}\in {\rm
int}S^{m-1}\setminus \{x^*\}$ the trajectory $\{x^{(n)}\}$ does not
converge.} \vskip 0.2 truecm

A skew-symmetric matrix $A$ is called transversal if all even order
leading (principal) minors are nonzero.

A Volterra QSO $V$ is called transversal if the corresponding
skew-symmetric matrix $A$ is transversal [14],[18],[36].

\vskip 0.2 truecm {\bf Problem 2.} {\sl Define concept of
transversality for arbitrary QSO and find necessary and sufficient
condition on matrix $\textbf{P}=\left(P_{ij,k}\right)$ of a QSO
under which the QSO is a transversal.}

\vskip 0.2 truecm Note that if a Volterra QSO is transversal then
the set $X=\{x\in S^{m-1}: V(x)=x\}$ of fixed points is a finite set
[14].

Let $U\equiv U_X$ be a neighborhood of the set $X$ and $\{x^{(n)}\}$
be an arbitrary trajectory. Denote
$$n_U=\left|\{j=1,...,n: x^{(j)}\in U\}\right|,$$
where $|M|$ denotes the number of elements in $M$.

Then it is known that
$$\lim_{n\to\infty}\frac{n_U}{n}=1$$
i.e. the trajectory a despondent part of the time will stay in the
neighborhood of the fixed points.

Denote $U=U_1\cup U_2\cup ... \cup U_t$, where $U_i$, $i=1,...,t$ is
the neighborhood of the fixed point $x_i$.

Thus the trajectory firstly visits the neighborhood of a fixed point
$x_{n_1}$ then it visits the neighborhood of a fixed point $x_{n_2}$
and so on.

The sequence $n_1, n_2, ... $ is called the itinerary (route-march)
of the trajectory $x^{(0)}, x^{(1)},...$ Since the set of fixed
points is a finite set, the numbers $n_1, n_2,...$ will repeat.
\vskip 0.2 truecm

{\bf Problem 3.} {\sl Is there a trajectory with periodic
itinerary?}\vskip 0.2 truecm

On the basis of numerical calculations Ulam [44] conjectured that
ergodic theorem holds for any QSO $V$, that is, the limit
$\lim_{n\to\infty}\sum_{k=0}^nV^k(x)$ exists for any $x\in S^{m-1}$.
In [45] Zakharevich proved that this conjecture is false in general.
In [3] the authors established a necessary condition for the ergodic
theorem to hold for the following class of Volterra QSOs $V:S^2\to
S^2$

$$\begin{array}{lll}
    x'=x\left(1+ay-bz\right), \\
   y'=y\left(1-ax+cz\right), \\
   z'=z\left(1+bx-cy\right),\\
\end{array}\eqno(9)
$$
where $a,b,c\in [-1,1].$

Note that if $a=b=c=1$ the QSO (9) coincides with an example
considered in [45]. \vskip 0.2 truecm

 {\bf Theorem 3.} [3] {\it If
the parameters $a,b,c$ for the Volterra QSO (9) have the same sign
and each is non-zero, then the ergodic theorem will fail for this
operator.} \vskip 0.2 truecm

 {\bf Problem 4.} {\sl Is the condition
of Theorem 3 sufficient for the ergodic theorem to hold?} \vskip 0.2
truecm

 {\bf Problem 5.} {\sl Find necessary and sufficient
conditions on matrix $A$ of Volterra QSO under which the ergodic
theorem is true for the Volterra QSO on $S^{m-1},\ \ m\geq 2$.}
\vskip 0.2 truecm
\section{The permuted Volterra QSO} Let $\tau$ be a cyclic
permutation on the set of indices $1,2,...,m$ and let $V$ be a
Volterra QSO. Define QSO $V_\tau$ by
$$V_\tau: x'_{\tau(j)}=x_j\left(1+\sum^m_{k=1}a_{jk}x_k\right), \
j=1,...,m, \eqno(10)$$ where $a_{jk}$ is defined in (8) (see
[12],[14],[16],[19]).

Note that QSO $V_\tau$ is a non-Volterra QSO iff $\tau\ne$id. \vskip
0.2 truecm

 {\bf Theorem 4.} [14] {\it For any automorphism
$W:S^{m-1}\to S^{m-1}$ there exists a permutation $\tau$ and a
Volterra QSO $V$ such that $W=V_\tau$.} \vskip 0.2 truecm

{\bf Corollary 1.} {\it The set of all quadratic automorphisms of
the simplex $S^{m-1}$ can be geometrically presented as the union of
$m!$ nonintersecting cubes of dimention $\frac{m(m-1)}{2}$.} \vskip
0.2 truecm

In [40] the behavior of trajectories of a non-Volterra automorphism
$V:S^2\to S^2$ are studied. \vskip 0.2 truecm

 {\bf Problem 6.} {\sl  Investigate
the asymptotic behavior of the trajectories of the operators
$V_\tau$ (automorphisms) for an arbitrary permutation $\tau$.}
\vskip 0.2 truecm

\section{$\ell$-Volterra QSO} Fix $\ell \in E$ and assume that
elements $P_{ij,k}$ of the matrix $\textbf{P}$ satisfy
$$ P_{ij,k}= 0 \ \ \mbox{if} \ \ k \not\in
\{i,j\}\ \ \ \mbox{for}\ \ \mbox{any}\ \ k\in \{1,...,\ell\},\ \ i,j
\in E;\eqno(11)$$
$$ P_{ij,k}> 0 \ \ \mbox{for at least one pair} \ \ (i,j),\ \
i\neq k,\ \ j\neq k \ \ \mbox{if} \ \  k \in
\{\ell+1,...,m\}.\eqno(11a)$$ \vskip 0.2 truecm

\textbf{Definition 2.} For any fixed $\ell \in E$, the QSO defined
by (4), (5), (11) and (11a) is called $\ell$-Volterra QSO. \vskip
0.2 truecm

 Denote by $\mathcal{V}_{\ell}$ the set of all
$\ell$-Volterra QSOs. \vskip 0.2 truecm

{\it Remarks.} 1. The condition (11a) guarantees that $\mathcal{V}
_{\ell_{1}} \bigcap \mathcal{V}_{\ell_{2}} = \emptyset $ for any
$\ell_{1}\neq \ell _{2}$.

2. Note that $\ell$-Volterra QSO is Volterra if and only if $\ell=
m$.

3. By Theorem 2 we know that there is no a periodic trajectory for
Volterra QSO. But for $\ell$-Volterra QSO there is such trajectories
(see Proposition 1 below). \vskip 0.2 truecm

 Let
$e_i=\left(\delta_{1i}, \delta_{2i},...,\delta_{mi}\right)\in
S^{m-1}$, $i=1,...,m$  be the vertices of $S^{m-1}$, where
$\delta_{ij}$ the Kronecker delta. \vskip 0.2 truecm

{\bf Proposition 1.}[41] 1) {\it For any set
$I_s=\{e_{i_1},...,e_{i_s}\}\subset \{e_{\ell+1},...,e_m\}, \ s\leq
m$, there exists a family $\mathcal{V} _{\ell}(I_s) \subset
\mathcal{V}_{\ell}$ such that $I_s$ is an $s$-cycle for every $V\in
\mathcal{V} _{\ell}(I_s).$}

2) For any $I_1,...,I_q\subset \{\ell+1,...,m\}$ such that $I_i\cap
I_j=\emptyset $ $(i\ne j, i,j=1,...,q)$, there exists a family
$\mathcal{V} _{\ell}(I_1,...,I_q) \subset \mathcal{V}_{\ell}$ such
that $\{e_i, i\in I_j\}$ $(j=1,...,q)$ is a $|I_j|$-cycle for every
$V\in \mathcal{V} _{\ell}(I_1,...,I_s)$. \vskip 0.2 truecm

{\bf Problem 7.} {\sl Find the set of all periodic trajectories of a
given $\ell-$Volterra QSO.} \vskip 0.2 truecm

In paper [41] the trajectories of an 1-Volterra and 2-Volterra QSOs
are studied. \vskip 0.2 truecm

 {\bf Problem 8.} {\sl Develop theory
of dynamical systems generated by a $\ell$-Volterra QSO. Find its
Lyapunov functions, the set of limit points of its trajectories
etc.} \vskip 0.2 truecm

 Note that in [4] a quasi-Volterra QSO was
considered, such a QSO is a particular case of $\ell$-Volterra QSO.

\section{Non-Volterra QSO as combination of a Volterra and a
non-Volterra operators}

In [15] it was considered the following family of QSOs $V_\lambda:
S^2\to S^2$: $V_\lambda=\lambda V_0+(1-\lambda)V_1, 0\leq\lambda\leq
1$, where $V_0(x)=\left(x_1^2+2x_1x_2, x_2^2+2x_2x_3,
x_3^2+2x_1x_3\right)$ is a Volterra QSO and
$V_1(x)=\left(x_1^2+2x_2x_3, x_2^2+2x_1x_3, x_3^2+2x_1x_2\right)$ is
a non-Volterra one.

Note that behavior of the trajectories of $V_0$ is very irregular
(see [25], [45]). It has fixed points $M_0=({1\over 3},{1\over
3},{1\over 3})$, $e_1, e_2, e_3$. The point $M_0$ is a repelling and
$e_i, i=1,2,3$ are saddle points. These four points, are also fixed
points for $V_1$ but $M_0$ is an attracting point for $V_1$. Thus
properties of $V_\lambda$ change depending on the parameter
$\lambda$. In [15] some examples of invariant curves and the set of
limit points of the trajectories of $V_\lambda$ are given. \vskip
0.2 truecm

 {\bf Problem 9.} {\sl For arbitrary two QSOs $V_1$ and
$V_2$ connect the properties of $V_\lambda=\lambda
V_1+(1-\lambda)V_2$, $\lambda\in [0,1]$ with properties of $V_1$ and
$V_2$.}

\section{F-QSO} Consider $E_0=E\cup \{0\}=\{0,1,...,m\}.$ Fix a set
$F\subset E$ and call this set the set of "females" and the set
$M=E\setminus F$ is called the set of "males".  The element $0$ will
play the role of an "empty-body".

Coefficients $P_{ij,k}$ of the matrix ${\mathbf P}$ we define as
follows
$$P_{ij,k}=\left\{\begin{array}{lll}
1, \ \ {\rm if} \ \ k=0, i,j\in F\cup \{0\} \ \ {\rm or} \ \ i,j\in
M\cup \{0\};\\
0, \ \ {\rm if} \ \ k\ne 0, i,j\in F\cup \{0\} \ \ {\rm or} \ \
i,j\in
M\cup \{0\};\\
\geq 0, \ \ {\rm if } \ \ i\in F, j\in M, \forall k.
\end{array}\right.\eqno (12)$$

The biological interpretation of the coefficients (12) is obvious:
the "child" $k$ can be born only if its parents are taken from
different classes $F$ and $M$. Generally, $P_{ij,0}$ can be strictly
positive for $i\in F$ and $j\in M$, which corresponds, for example,
to the case in which "female" $i$ with "male" $j$ cannot have a
"child", because one of them is ill or both are.

\vskip 0.3 truecm

{\bf Definition 3.} {\it  For any fixed $F\subset E$, the QSO
defined by (4),(5) and (12) is called the $F-$ quadratic stochastic
operator ($F$-QSO).}

\vskip 0.3 truecm

{\it Remarks.} 1. Any $F-$ QSO is non- Volterra, because
$P_{ii,0}=1$ for any $i\ne 0.$

2. For $m=1$ there is a unique $F-$QSO (independently of $F=\{1\}$
and $F=\emptyset$) which is constant i.e., $V(x)=(1,0)$ for any
$x\in S^1.$ \vskip 0.2 truecm

{\bf Theorem 5.} [39] {\it Any $F$-QSO has a unique fixed point
$(1,0,...,0)$ (with $m$ zeros). Besides, for any $x^0\in S^m$, the
trajectory $\{x^{(n)}\}$ tends to this fixed point exponentially
rapidly.} \vskip 0.2 truecm

{\bf Problem 10.} {\sl Consider a partition $\xi=\{E_1,...,E_q\}$ of
$E$ i.e., $E=E_1\cup...\cup E_q$, $E_i\cap E_j=\emptyset$, $i\ne j.$
Assume $P_{ij,k}=0$ if $i,j\in E_p$, for $p=1,...,q$. Call the
corresponding operator a $\xi-$QSO. Is an analogue of Theorem 5 true
for any $\xi-$QSO?}

\section{Strictly non-Volterra QSO} Recently in [40] a new class of
non-Volterra QSOs have been introduced. These QSOs satisfy
$$P_{ij,k}=0 \ \ \mbox{if} \ \ k\in\{i,j\}, \ \forall i,j,k\in E.\eqno(13)$$
Such an operator is called strictly non-Volterra QSO. One can easily
check that the strictly non-Volterra operators exist only for $m\geq
3$.

An arbitrary strictly non-Volterra QSO defined on $S^2$ (i.e.,
$m=3$) has the form:

$$\begin{array}{lll}
    x'=\alpha y^2+cz^2+2yz, \\
   y'=ax^2+dz^2+2xz, \\
   z'=bx^2+\beta y^2+2xy,\\
\end{array}\eqno(14)
$$
where $$a,b,c,d, \alpha,\beta\geq 0, \ \
a+b=c+d=\alpha+\beta=1.\eqno(15)$$

\vskip 0.2 truecm

{\bf Theorem 6.} [40] 1) {\it For any values of parameters
$a,b,c,d,\alpha, \beta$ with (15) the operator (14) has a unique
fixed point. Moreover the fixed point is not attractive.}

2) {\it The QSO (14) has 2-cycles and 3-cycles depending on the
parameters (15).}

\vskip 0.2 truecm
 {\bf Problem 11.} {\sl Is Theorem 6 true for
$m\geq 4$?}

\section{Regularity of QSO} In [17] the authors consider an
arbitrary QSO $V:S^{m-1}\to S^{m-1}$ with matrix
$\textbf{P}=(P_{ij,k})$ and studied the problem of finding the
smallest $\alpha_m$ such that the condition $P_{ij,k}>\alpha_m$
implies the regularity of $V$.

\vskip 0.2 truecm
{\bf Theorem 7.} [17] 1) {\it If $P_{ij,k}>{1\over
2m}$ then $V$ is regular.}

2) $\alpha_2={1\over 2}(3-\sqrt{7}).$

\vskip 0.2 truecm
 {\bf Problem 12.} {\sl Find exact values of
$\alpha_m$ for any $m\geq 3$.}

\section{Quadratic bistochastic operators} Let $x\in S^{m-1}$.
Denote by $x_\downarrow$ the point
$x_\downarrow=(x_{[1]},...,x_{[m]})\in S^{m-1}$, where
$x_{[1]}\geq...\geq x_{[m]}$ are the coordinates of $x$ in
non-increasing order.

If $x,y\in S^{m-1}$ and
$$\sum_{i=1}^kx_{[i]}\leq \sum^k_{i=1}y_{[i]}, \ k=1,...,m,$$
then we say that $y$ majorizes $x$ and write $x\prec y$.

As is known [27], $x\prec y$ iff there is a doubly stochastic
(bistochastic) matrix $B$ such that $x=By$. Therefore, if $B$ is a
bistochastic matrix, then $Bx\prec x$ for any point $x\in S^{m-1}$.

In [21] it was considered more general definition: \vskip 0.2 truecm

{\bf Definition 4.} An arbitrary continuous operator $V:S^{m-1}\to
S^{m-1}$ satisfying the condition
$$V(x)\prec x, \ \ x\in S^{m-1}\eqno(16)$$
is called a bistochastic operator. \vskip 0.2 truecm

{\bf Theorem 8.} [21],[22] 1) {\it If $V:S^{m-1}\to S^{m-1}$ is a
bistochastic operator, then the coefficients $P_{ij,k}$ satisfy the
conditions}
$$\sum_{i,j=1}^mP_{ij,k}=m, \ \ \forall k=1,...,m;\eqno(a)$$
$$\sum_{j=1}^mP_{ij,k}\geq {1\over 2},\ \ \forall i,k=1,...,m;\eqno(b)$$
$$\sum_{i,j\in I}P_{ij,k}\leq t, \ \ \forall t,k=1,...,m,\eqno(c)$$
{\it where $I=\{i_1,...,i_t\}$ is an arbitrary subset of
$\{1,..,m\}$ containing $t$ elements.}

2) {\it If (c) holds for a QSO $V$ then it is a bistochastic.}

\vskip 0.2 truecm

 Let $\textbf{B}$ be the set of all bistochastic
quadratic operators acting in $S^{m-1}$. The set $\textbf{B}$ can be
regarded as a polyhedron in an ${m(m^2-1)\over 2}$-dimensional
space. Let Extr$\textbf{B}$ be the set of extreme points of
$\textbf{B}$. \vskip 0.2 truecm

 {\bf Theorem 9.}[20],[21] {\it If
$V\in$Extr$\textbf{B}$, then}
$$P_{ii,k}=0\ \ \mbox{or}\ \ 1; \eqno(d)$$
$$P_{ij,k}=0, {1\over 2}\ \ \mbox{or}\ \ 1, \ \ \mbox{for} \ \ i\ne j. \eqno(d)$$

\vskip 0.2 truecm
 Note that the converse assertion of the Theorem is
false [21].

In [22] an  analogue of Birkhoff's theorem is proved. \vskip 0.2
truecm
 {\bf Problem 13.} {\sl Investigate the behavior of
trajectories of the bistochastic quadratic operators.}

\section{Surjective QSOs} In [5] and [28] a description of
surjective QSOs defined on $S^{m-1}$ for $m=2,3,4$  and
classification of extreme points of the set of such operators are
given. \vskip 0.2 truecm

 {\bf Problem 14.} {\sl Describe the set of
all surjective QSOs defined on $S^{m-1}$ for any $m\geq 5$.}

\section{Construction of QSO}  In papers [6],[7]
a constructive description of $\textbf{P}$ (i.e. QSO) is given. The
construction depends on cardinality of $E$, namely two cases: (i)
$E$ is finite (ii) $E$ is a continual set were separately
considered. Note that for the second case one of the key problem is
to determine the set of coefficients of heredity which is already
infinite dimensional; the second problem is to investigate the
quadratic operator which corresponds to this set of coefficients. By
the construction the operator $V$ depends on a probability measure
$\mu$ being defined on a measurable space $(E,{\cal F}).$

Recall the construction for finite $E=\{1,...,m\}$.

Let $G=(\L, L)$ be a finite graph without loops and multiple edges,
where $\L$ is the set of vertexes and $L$ is the set of edges of the
graph.

Furthermore, let $\Phi$ be a finite set, called the set of alleles
(in problems of statistical mechanics, $\P$ is called the  range of
spin). The function $\s:\L\to\P$ is called a cell (in mechanics it
is called configuration). Denote by $\Omega$ the set of all cells,
this set corresponds to $E$. Let $S(\L, \P)$ be the set of all
probability measures defined on the finite set $\Omega.$

Let $\{\L_i, i=1,...,q\}$ be the set of maximal connected subgraphs
(components) of the graph $G.$ For any $M\subset \L$ and $\s\in
\Omega$ denote $\s(M)=\{\s(x): x\in M\}.$ Fix two cells $\s_1,
\s_2\in \Omega,$ and put
$$\Omega(G, \s_1,\s_2)=\{\s\in \Omega: \s(\L_i)=\s_1(\L_i) \ \ {\rm or} \ \
\s(\L_i)=\s_2(\L_i) \ \ {\rm for \ \ all} \ \ i=1,...,m\}.$$ Now let
$\mu\in S(\L,\P)$ be a probability measure defined on $\Omega$ such
that $\mu(\s)>0$ for any cell $\s\in \Omega;$ i.e $\mu$ is a Gibbs
measure with some potential [37]. The heredity coefficients
$P_{\s_1\s_2,\s}$ are defined as
$$P_{\s_1\s_2,\s}=\left\{\begin{array}{ll}
{\mu(\s)\over \mu(\Omega(G,\s_1,\s_2))}, \ \ {\rm if} \ \ \s\in \Omega(G,\s_1,\s_2),\\
0 \ \ {\rm otherwise}.\\
\end{array}\right.\eqno (17)$$
Obviously, $P_{\s_1\s_2,\s}\geq 0,$
$P_{\s_1\s_2,\s}=P_{\s_2\s_1,\s}$ and $\sum_{\s\in \Omega}
P_{\s_1\s_2,\s}=1$ for all $\s_1,\s_2\in \Omega.$

The QSO $V\equiv V_\mu$ acting on the simplex $S(\L,\P)$ and
determined by coefficients (17) is defined as follows: for an
arbitrary measure $\l\in S(\L,\P)$, the measure $V(\l)=\l'\in
S(\L,\P)$ is defined by the equality
$$ \l'(\s)=
\sum_{\s_1,\s_2\in \Omega} P_{\s_1\s_2,\s} \l(\s_1)\l(\s_2)
\eqno(18)$$ for any cell $\s\in \Omega.$

 \vskip 0.3 truecm
{\bf Theorem 10.} [6] {\it The QSO (18) is Volterra if and only if
the graph $G$ is connected.} \vskip 0.3 truecm

Thus if $\Phi$, $G$ and $\mu$ are given then we can constuct a QSO
corresponding to these objects. In [6],[29] several examples of
$\Phi$, $G$ and $\mu$ are considered and the trajectories of
corresponding QSOs are studied. \vskip 0.2 truecm

Note that the construction above does not give all possible QSOs. So
the following problem is interesting. \vskip 0.2 truecm

 {\bf Problem 15.} {\sl Describe the class of QSOs which can be obtained by the
construction.}

\vskip 0.2 truecm In [7] also constructively described QSOs which
act to the set of all probability measures on some measurable space
$(E,{\cal F})$ where $E$ is a uncountable set. This construction
depends on a Gibbs measure $\mu$ (see [37]). The behavior of
trajectories of such operators were studied. These investigations
allows to a natural introduction of thermodynamics in studying some
models of heredity. More precisely, if $E$ is continual set then one
can associate the Gibbs measure $\mu$ by a Hamiltonian $H$ (defined
on $E$) and temperature $T>0$ [37]. It is known that depending on
the Hamiltonian and the values of the temperature the measure $\mu$
can be non unique. In this case there is a phase transition of the
physical system with the Hamiltonian $H$.

In [7] for $q$-state Potts Hamiltonian when the temperature is low
enough, it is proven that any trajectory of the QSO constructed by a
Gibbs measure $\mu_i$, $i=1,...,q$ of the Potts Hamiltonian  tends
to the measure $\mu_i$. In other words, any trajectory of the QSO
generated by a Gibbs measure of the Potts model converges to this
measure.
\vskip 0.2 truecm
{\bf Problem 16.} (by N.N.Ganikhodjaev)
{\sl How the thermodynamics (the phase transition) will effect to
behavior of the trajectories of a QSO corresponding to a Gibbs
measure of the Hamiltonian $H$?}

\section{Non-Volterra QSO generated by a product measure}

In [38] it was shown that if $\mu$ is  the product of probability
measures being defined on each maximal connected subgraphs of $G$
then corresponding non-Volterra operator can be reduced to $q$
number (where $q$ is the number of maximal connected subgraphs of
$G$) of Volterra operators defined on the maximal connected
subgraphs.

Let $G=(\L,L)$ be a finite graph and $\{\L_i, i=1,...,q\}$ the set
of all maximal connected subgraphs of $G$. Denote by
$\Omega_i=\P^{\L_i}$ the set of all configurations defined on
$\L_i,$ $i=1,...,q.$ Let $\mu_i$ be a probability measure defined on
$\Omega_i,$ such that $\mu_i(\s)>0$ for any $\s\in \Omega_i,$
$i=1,...,q.$

Consider probability measure $\mu$ on
$\Omega=\Omega_1\times\dots\times \Omega_q$ defined by

$$ \mu(\s)=\prod^q_{i=1}\mu_i(\s_i), \eqno (19)$$
where $\s=(\s_1,...,\s_q),$ with $\s_i\in \Omega_i, i=1,...,q.$

According to Theorem 10 if $q=1$ then QSO constructed on $G$ is
Volterra QSO. \vskip 0.4 truecm

{\bf Theorem 11.} [38] {\it The QSO constructed by the construction
(18) with respect to measure (19) is reducible to $q$ separate
Volterra QSOs.}

\vskip 0.2 truecm This result allows to study a wide class of
non-Volterra operators in the framework of the well known theory of
Volterra quadratic stochastic operators.

\vskip 0.2 truecm {\bf Problem 17.} {\sl Describe the set of all
non-Volterra QSOs which are reducible to several Volterra
QSOs.}\vskip 0.2 truecm

{\bf Problem 18.} {\sl Find a measure $\mu$ different from (19) such
that the non-Volterra QSO corresponding to $\mu$ can be investigated
in the framework of a well known theory of QSOs.}

\section{Trajectories with historic behavior} The problem which we
shall discuss here is a particular case of the problem stated in
[43].

Consider a QSO $V:S^{m-1}\to S^{m-1}$. We say that a trajectory
$\{x, V(x), V^2(x), ...\}$ has historic behavior if for some
continuous function $f:S^{m-1}\to R$ the average
$$\lim_{n\to \infty} {1\over n+1}\sum^n_{i=0}f(V^i(x))$$
does not exist.

If this limit does not exist, it follows that "partial averages"
${1\over n+1}\sum^n_{i=0}f(V^i(x))$ keep changing considerable so
that their values give information about the epoch to which $n$
belongs: they have a history [43]. \vskip 0.2 truecm

 {\bf Problem
19.} {\sl Find a class of QSOs such that the set of initial states
which give rise to trajectories with historic behavior has positive
Lebesque measure.} \vskip 0.2 truecm

 Similar problem was discussed
by Ruelle in [40].

\section{A generalization of Volterra QSO}
Consider QSO (4), (5) with additional condition
$$P_{ij,k}=a_{ik}b_{jk}, \ \ \forall i,j,k \in E\eqno(20)$$
where $a_{ik}, b_{jk}\in R$ entries of matrices $A=(a_{ik})$ and
$B=(b_{jk})$ such that conditions (5) are satisfied for the
coefficients (20).

Then the QSO $V$ corresponding to the coefficients (20) has the form
$$x'_k=(V(x))_k=(A(x))_k\cdot (B(x))_k, \eqno(21) $$
where $$(A(x))_k=\sum^m_{i=1}a_{ik}x_i, \ \
(B(x))_k=\sum^m_{j=1}b_{jk}x_j.$$
 Note that if $A$ (or $B$) is the identity matrix then operator (21)
is a Volterra QSO.

\vskip 0.2 truecm {\bf Problem 20.} {\sl Develop theory of QSOs
defined by (21).} \vskip 0.2 truecm

\section{Bernstein's problem} The Bernstein problem [25],[26] is related to
a fundamental statement of population genetics, the so-called
stationarity principle. This principle holds provided that the
Mendel law is assumed, but it is consistent with other mechanisms of
heredity. An adequate mathematical problem is as follows. QSO $V$ is
a Bernstein mapping if $V^2=V$. This property is just the
stationarity principle. This property also is known as
Hardy-Weinberg law [23]. The problem is to describe all Bernstein
mappings explicitly. The case $m\leq 2$ is mathematically trivial
and biologically not interesting. Bernstein [1] solved the above
problem for the case $n=3$ and obtained some results for $n\geq 4$.
In works by Lyubich (see e.g.[25],[26]) the Bernstein problem was
solved for all $m$ under the regularity assumption. The regularity
means that $V(x)$ depends only on the values $f(x)$, where $f$ runs
over all invariant linear forms. In investigations by Lyubich [25],
the algebra $A_V$ with the structure constants $P_{ij,k}$ played a
very important role. Since $V(x)=x^2$, the Bernstein property of $V$
is equivalent ti the identity
$$(x^2)^2=s^2(x)x^2.$$
This identity means that $A_V$ is a Bernstein algebra with respect
to the algebra homomoprphism $s:A_V\to R$. The mapping $V$ is
regular iff the identity
$$x^2y=s(x)xy$$
holds in the algebra $A_V$, by definition, this identity means that
$A_V$ is regular.
 \vskip 0.2 truecm
 {\bf Problem 21.} {\sl Describe all QSOs which satisfy
$V^r(x)=V(x)$ for any $x\in S^{m-1}$ and some $r\geq 2$.}

\section{Topological conjugacy}

{\bf Definition 5.} Let $V_1:S^{m-1}\to S^{m-1}$ and $V_2:S^{m-1}\to
S^{m-1}$ be two QSOs with coefficients $P^{(1)}_{ij,k}$ and
$P^{(2)}_{ij,k}$ respectively. $V_1$ and $V_2$ are said to be
topologically conjugate if there exists a homeomorphism
$h:S^{m-1}\to S^{m-1}$ such that, $h\circ V_1=V_2\circ h$. The
homeomorphism $h$ is called a topological conjugacy.

Mappings which are topologically conjugate are completely equivalent
in terms of their dynamics [2]. \vskip 0.2 truecm

{\bf Definition 6.} A polynomial $f(P_{ij,k})$ is called an
indicator if from the topologically conjugateness of $V_1$ and $V_2$
it follows that
$$\alpha_f\leq f(P^{(1)}_{ij,k})\leq \beta_f\ \ \mbox{and} \ \
\alpha_f\leq f(P^{(2)}_{ij,k})\leq \beta_f,$$ where $\alpha_f,
\beta_f\in R$. \vskip 0.2 truecm

 {\bf Definition 7.} A system
$f_1,...,f_t$ of indicators is called complete if from
$$\alpha_{f_n}\leq f_n(P^{(1)}_{ij,k})\leq \beta_{f_n}\ \ \mbox{and} \ \
\alpha_{f_n}\leq f_n(P^{(2)}_{ij,k})\leq \beta_{f_n},$$ for any
$n=1,...,t$ it follows the topologically conjugateness of $V_1$ and
$V_2$. \vskip 0.2 truecm

 A minimal complete system of indicators is
called a basis. \vskip 0.2 truecm

 {\bf Problem 22.} {\sl Does there
exist a finite complete system of indicators? Find the basis of the
system of indicators.} \vskip 0.2 truecm

 {\it Remark.} There are
many papers devoted to (quantum) quadratic stochastic processes
[8],[9],[31]-[33] and to an infinite dimensional Volterra quadratic
operators [34], [35]. Even exotic directions, such as the p-adic
QSOs has been considered in [24]. Many other problems also discussed
in [30]. \vskip 0.3 truecm

 {\bf Acknowledgments.} This work was done
within the scheme of Junior Associate at the ICTP, Trieste, Italy
and UAR thanks ICTP for providing financial support and all
facilities (in February - April 2009). \vskip 0.4 truecm

{\bf References}\footnote{Ganikhodjaev(Ganikhodzhaev) N.N. is Nasir
Nabievich and Ganikhodjaev(Ganikhodzhaev) R.N. is Rasul Nabievich.
Rasul is brother of Nasir.}

1. Bernstein S.N. The solution of a mathematical problem related to
the theory of heredity.{\it Uchn. Zapiski NI Kaf. Ukr. Otd. Mat.}
1924. no. 1., 83-115.(Russian)

2.  Devaney R. L., An introduction to chaotic dynamical system, {\it
Westview Press,} 2003.

3. Ganikhodzhaev, N. N.; Zanin, D. V. On a necessary condition for
the ergodicity of quadratic operators defined on a two-dimensional
simplex. {\it Russian Math. Surveys.} 59 (2004), no. 3, 571-572.

4. Ganikhodzhaev, N. N.; Mukhitdinov, R. T. On a class of
non-Volterra quadratic operators. {\it Uzbek. Mat. Zh.} (2003), no.
3-4, 9--12.(Russian)

5. Ganikhodjaev ,N.N., Mukhitdinov R.T. Extreme points of a set of
quadratic operators on the simplices $S^1$ and $S^2$. {\it Uzbek.
Mat.Zh.} (1999), no.3, 35-43.(Russian)

6. Ganikhodjaev N.N. An application of the theory of Gibbs
distributions to mathematical genetics. {\it Doklady Math}. 61
(2000), 321-323.

7. Ganikhodjaev, N. N.; Rozikov, U. A. On quadratic stochastic
operators generated by Gibbs distributions. {\it Regul. Chaotic
Dyn.} 11 (2006), no. 4, 467-473.

8. Ganikhodjaev, N. N.; Akin, H.; Mukhamedov, F.M. On the ergodic
principle for Markov and quadratic stochastic processes and its
relations. {\it Linear Algebra Appl.} 416 (2006), no. 2-3, 730-741.

9. Ganikhodzhaev, N. N.; Mukhamedov, F. M. On the ergodic properties
of discrete quadratic stochastic processes defined on von Neumann
algebras. {\it Izvestiya: Math.} 64 (2000), no. 5, 873-890.

10. Ganikhodzhaev R. N. Quadratic stochastic operators, Lyapunov
functions and tournaments.{\it Acad. Sci. Sb. Math.} 76 (1993), no.
2, 489-506.

11. Ganikhodzhaev, R. N. A chart of fixed points and Lyapunov
functions for a class of discrete dynamical systems. {\it Math.
Notes.} 56 (1994), no. 5-6, 1125-1131.

12. Ganikhodzhaev, R. N.; Karimov, A. Z. Mappings generated by a
cyclic permutation of the components of Volterra quadratic
stochastic operators whose coefficients are equal in absolute
magnitude. {\it Uzbek. Mat. Zh.} No. 4 (2000), 16-21.(Russian)

13. Ganikhodzhaev, R. N.; Eshmamatova, D. B. On the structure and
properties of charts of fixed points of quadratic stochastic
operators of Volterra type. {\it Uzbek. Mat. Zh.} No. 5-6 (2000),
7-11. (Russian)

14. Ganikhodzhaev, R. N.; Eshmamatova, D. B. Quadratic automorphisms
of a simplex and the asymptotic behavior of their trajectories. {\it
Vladikavkaz. Mat. Zh.} 8 (2006), no. 2, 12-28.(Russian)

15. Ganikhodzhaev, R. N. A family of quadratic stochastic operators
that act in $S\sp 2$. {\it Dokl. Akad. Nauk UzSSR}. (1989), no. 1,
3-5.(Russian)

16. Ganikhodzhaev, R. N.; Dzhurabaev, A. M. The set of equilibrium
states of quadratic stochastic operators of type $V\sb {\pi}$. {\it
Uzbek. Mat. Zh.}  No. 3 (1998), 23-27.(Russian)

17. Ganikhodzhaev, R. N.; Sarymsakov, A. T. A simple criterion for
regularity of quadratic stochastic operators. {\it Dokl. Akad. Nauk
UzSSR.} (1988), no. 11, 5-6.(Russian)

18. Ganikhodzhaev, R. N.; Abdirakhmanova, R. E. Description of
quadratic automorphisms of a finite-dimensional simplex. {\it Uzbek.
Mat. Zh. No.} 1 (2002), 7-16.(Russian)

19. Ganikhodzhaev, R. N.; Abdirakhmanova, R. E. Fixed and periodic
points of quadratic automorphisms of non-Volterra type. {\it Uzbek.
Mat. Zh.} (2002), no. 2, 6-13.(Russian)

20. Ganikhodzhaev, R. N.; Eshniyazov, A. I. Bistochastic quadratic
operators. {\it Uzbek. Mat. Zh.} (2004), no. 3, 29-34.(Russian)

21. Ganikhodzhaev, R. N. On the definition of quadratic bistochastic
operators. {\it Russian Math. Surveys} 48 (1993), no. 4, 244-246.

22. Ganikhodzhaev, R. N.; Shahidi F. Doubly stochastic quadratic
operators and Birkhoff's problem. {Arxiv:} 0802.1100 [math.FA] 2008.

23. Hofbaver J., Sigmund K., The theory of evolution and dynamical
systems, \textit{Cambridge Univ. Press}, 1988.

24. Khamraev, M. M. On $p$-adic dynamical systems associated with
Volterra type quadratic operators of dimension 2. {\it Uzbek. Mat.
Zh.} (2005), no. 1, 88-96.

25. Lyubich Yu. I., Mathematical structures in population genetics,
{\sl Biomathematics}, Springer-Verlag, {\bf 22} (1992).

26. Lyubich Yu.I., Ultranormal case of the Bernstein problem. {\it
Func. Anal. Appl.} 31 (1997), no.1, 60-62.

27. Marshall, A.W, Olkin I. Inequalities: theory of majorization and
its applications, {\it Acad. Press.} New York. 1979.

28. Meyliev, Kh. Zh. Description of surjective quadratic operators
and classification of the extreme points of a set of quadratic
operators defined on $S\sp 3$. {\it Uzbek. Mat. Zh.} (1997), no.3,
39-48.(Russian)

29. Meyliev, Kh. Zh.; Mukhitdinov, R. T.; Rozikov, U. A. On two
classes of quadratic operators that correspond to Potts models and
$\lambda$-models. {\it Uzbek. Mat. Zh. No.} (2001), no.1,
23-28.(Russian)

30. F.Mosconi, et al. Some nonlinear challenges in biology. {\it
Nonlinearity.} 21 (2008) T131-T147.

31. Mukhamedov, F. M. On the decomposition of quantum quadratic
stochastic processes into layer-Markov processes defined on von
Neumann algebras. {\it Izvestiya: Math.} 68 (2004), no. 5,
1009-1024.

32.  Mukhamedov, F. M. On a regularity condition for quantum
quadratic stochastic processes. {\it Ukrainian Math. J.} 53 (2001),
no. 10, 1657-1672.

33. Mukhamedov, F. M. Ergodic properties of quadratic discrete
dynamical systems on $C\sp *$-algebras. {\it Methods Funct. Anal.
Topology} 7 (2001), no. 1, 63-75.

34. Mukhamedov, Farruh; Akin, Hasan; Temir, Seyit On infinite
dimensional quadratic Volterra operators. {\it J. Math. Anal. Appl.}
310 (2005), no. 2, 533-556.

35. Mukhamedov, F. M. On infinite-dimensional quadratic Volterra
operators. {\it Russian Math. Surveys} 55 (2000), no. 6, 1161-1162.

36. Mukhamedov, F.M.; Saburov M. On homotopy of Volterrian quadratic
stochastic operators. Arxiv:0712.2891 [math.DS], 2007.

37. Preston C. {\it Gibbs measures on countable sets.} Cambridge
Univ. Press.,London 1974.

38. Rozikov U.A., Shamsiddinov N.B. On non-Volterra quadratic
stochastic operators generated by a product measure. {\it Stochastic
Anal. Appl.} 27 (2009), no.2, p.353-362.

39. Rozikov U.A., Zhamilov U.U. On F-quadratic stochastic operators.
{\it Math. Notes.} 83 (2008), no.4, 554-559.

40. Rozikov U.A., Zhamilov U.U. On dynamics of strictly non-Volterra
quadratic operators defined on the two dimensional simplex. To
appear in {\it Sbornik: Math.}

41. Rozikov U.A., Zada A. On $\ell$-Volterra quadratic stochastic
operators. {\it Doklady Math.} 79 (2009) no.1, 32-34.

42. Ruelle D. Historic behavior in smoth dynamical systems. {\it
Global Anal. Dynamical Syst.} ed. H.W.Broer et al. 2001. (Bristol:
Institute of Physics Publishing)

43. Takens F. Orbits with historic behavior, or non-existence of
averages. {\it Nonlinearity.} 21 (2008), T33-T36.

44. Ulam, S.M. {\it Problems in Modern Math.,} New York; Wiley,
1964.

45. Zakharevich, M.I. the behavior of trajectories and the ergodic
hypothesis for quadratic mappings of a simplex. {\it Russian Math.
Surveys.} 33 (1978), 207-208.

\end{document}